\documentclass[11pt]{article}
\usepackage[numbers,sort&compress]{natbib}
\usepackage[colorlinks,
            linkcolor=blue,
            anchorcolor=green,
            citecolor=magenta
            ]{hyperref}
\usepackage{mathrsfs,amssymb,amsfonts}
\usepackage{graphicx}
\usepackage{amsmath,amssymb,latexsym,color}
\usepackage[mathscr]{eucal}
\usepackage[numbers]{natbib}
\usepackage{dashrule}
\usepackage{pifont}

\def\thefootnote{\fnsymbol{footnote}}
\newtheorem{thm}{Theorem}[section]

\newtheorem{prop}[thm]{Proposition}
\newtheorem{lem}[thm]{Lemma}

\newtheorem{cor}[thm]{Corollary}
\newtheorem{example}{Example}
\newtheorem{remark}{Remark}[section]

\newcommand{\gc}[2]{{{#1}\brack #2}}

\def\msf{\,\mathscr{F}}
\def\msm{\,\mathscr{M}}

\def\msa{\,\mathscr{A}}
\def\msb{\,\mathscr{B}}
\def\msc{\,\mathscr{C}}

\newcommand{\proof}{{\it Proof.\quad}}
\newcommand{\qed}{\hfill\Box\medskip}
\begin{document}

\renewcommand{\abovewithdelims}[2]{
\genfrac{[}{]}{0pt}{}{#1}{#2}}
\renewcommand{\thefootnote}{}

\title{\bf Intersecting families of vector spaces with maximum covering number}

\author{Chao Gong\quad Benjian Lv\quad Kaishun Wang\\
{\footnotesize   \em  Sch. Math. Sci. {\rm \&} Lab. Math. Com. Sys., Beijing Normal University, Beijing, 100875,  China}
}
 \date{}
 \date{}
 \maketitle

\begin{abstract}
Let $V$ be an $n$-dimensional vector space over the finite field $\mathbb{F}_q$. Suppose that $\mathscr{F}$ is an intersecting family of $m$-dimensional subspaces of $V$. The covering number of $\msf$ is the minimum dimension of a subspace of $V$ which intersects all elements of $\mathscr{F}$. In this paper, we give the tight upper bound for the size of $\mathscr{F}$ whose covering number is $m$, and describe the structure of $\mathscr{F}$ which reaches the upper bound. Moreover, we determine the structure of an maximum intersecting family of singular linear space with the maximum covering number. 

\medskip



\footnote{ E-mail address:\\
gongchaomathc@163.com (C. Gong), bjlv@bnu.edu.cn (B. Lv),
wangks@bnu.edu.cn (K. Wang).}

\end{abstract}
\section{Introduction}
Let $X$ be an $n$-element set and ${X\choose m}$ denote the family of all $m$-subsets of $X$. A family $\msf\subseteq\binom{X}{m}$ is called \emph{intersecting} if for all $F_{1},F_{2}\in\msf$ we have $F_{1}\cap F_{2}\ne \emptyset$. The \emph{covering number} $\tau(\msf)$ is the minimum size of a set that meets all $F\in\msf$. We say that $\msf$ is \emph{trivial} if $\tau(\msf)=1$. Erd\H{o}s, Ko and Rado \cite{KO} determined the maximum size of an intersecting family $\msf$ with $n>2m$ and showed that any intersecting family with maximum size is trivial. In 1967, Hilton and Milner \cite{AE} determined the maximum size of a non-trivial intersecting family $\msf$ with $n>2m$. In 1986, Frankl and F\"{u}redi \cite{FF} gave a new proof using the shifting technique. Their results showed that any non-trivial intersecting family with maximum size must have covering number $2$. Over the years, there have been some results about covering number of intersecting family, see \cite{FRANKL, SM, MM, PKN, PFRANKL}. 

For any intersecting family $\msf\subseteq{X\choose m}$, we have $1\leq \tau(\msf)\leq m$. Let
$$r(m)=\max\left\{|\msf|: \msf\subseteq\binom{X}{m}~\text{is an intersecting family with}~\tau(\msf)=m\right\}.$$
The bounds for $r(m)$ were obtained in \cite{AT,PK}.
\begin{thm}{\rm(\cite{PK, AT})}
$\left(\frac{m}{2}\right)^{m-1}\leq r(m)\leq (1+O(1))m^{m-1}$.
\end{thm}

Let $V$ be an $n$-dimensional vector space over the finite field $\mathbb{F}_q$ and $\gc{V}{m}$ denote the family of all $m$-subspaces of $V$. For $n,m\in\mathbb{Z}^+$, define the \emph{Gaussian binomial coefficient} by
$$\gc{n}{m}:=\prod_{0\leq i<m}\frac{q^{n-i}-1}{q^{m-i}-1}.$$
Note that the size of $\gc{V}{m}$ is $\gc{n}{m}$.

For two subspaces $A,B\subseteq V$, we say that $A$ \emph{intersects} $B$ if $\dim(A\cap B)\geq 1$. A family $\msf\subseteq\gc{V}{m}$ is called \emph{intersecting} if $A$ intersects $B$ for all $A,B\in\msf$.
For any $\msf\subseteq\gc{V}{m}$, the \emph{covering number} $\tau(\msf)$ is the minimum dimension of a subspace of $V$ that intersects all elements of $\msf$. For any intersecting family $\msf\subseteq\gc{V}{m}$, we have $1\leq \tau(\msf)\leq m$. We say that $\msf$ is \emph{trivial} if $\tau(\msf)=1$. The maximum size of an intersecting family $\msf$ was determined in \cite{PL, PR, Hsieh, WN} using different techniques. It was showed that $\msf$ with maximum size has covering number $1$. Blokhuis et al. \cite{AB} determined the maximum size of intersecting family whose covering number is more than $1$. 

In this paper, we determine the maximum size of an intersecting family $\msf\subseteq\gc{V}{m}$ with $\tau(\msf)=m$. Our main result is as follow.
\begin{thm}\label{1.2}
Suppose $q\geq m\geq 2$ and $n\geq 2m-1$. Let $V$ be an $n$-dimensional vector space over $\mathbb{F}_q$ and $\msf\subseteq\gc{V}{m}$ be an intersecting family with $\tau(\msf)=m$. Then $|\msf|\leq\gc{2m-1}{m}$ and equality holds if and only if $\msf=\gc{X}{m}$ for some $X\in\gc{V}{2m-1}$.
\end{thm}

In Section 2, we give a proof of Theorem \ref{1.2}. In Section 3, we determine the maximum size of an intersecting family of singular linear spaces with maximum covering number. 

\section{Proof of Theorem \ref{1.2}}
Let $V$ be an $n$-dimensional vector space over $\mathbb{F}_q$ and $\msf\subseteq\gc{V}{m}$ be an intersecting family with $\tau(\msf)=m$. For any subspaces $A$ and $B$ of $V$, write $\msf_A=\{F\in\msf:~A\subseteq F\}$, and write $A\cap B=0$ if $\dim(A\cap B)=0$.

As is pointed out in \cite{AB}, Theorem \ref{1.2} holds for $m=2$. In the following, we always assume that $m\geq 3$. 

\begin{lem}{\rm(\cite[Lemma 2.4]{AB})}\label{2.1}
If an $s$-subspace $S$ does not intersect each element of $\msf$, then there exists an $(s+1)$-subspace $T$ such that $|\msf_T|\geq|\msf_S|/\gc{m}{1}$.
\end{lem}

\begin{cor}\label{2.2}
Suppose $s\leq t\leq m$. For any $S\in\gc{V}{s}$, there exists $T\in\gc{V}{t}$ such that $|\msf_S|\leq\gc{m}{1}^{t-s}|\msf_T|$. Moreover, we have $|\msf_S|\leq\gc{m}{1}^{m-s}$.
\end{cor}
\proof
If $t=s$, then pick $T=S$. Suppose $s<t$. For any $I\in\gc{V}{i}$ with $s\leq i<m$, there exists an $F\in\msf$ such that $F\cap I=0$. Applying Lemma \ref{2.1} $t-s$ times, we get the former inequality. Since $|\msf_W|\leq 1$ for any $W\subseteq\gc{V}{m}$, the later inequality holds.
$\qed$

Denote $X=\sum_{F\in\msf}F$. For any $R\in\gc{X}{m-1}$, since $\tau(\msf)=m$, there exists an $F\in\msf$ such that $F\cap R=0$, which implies that
\begin{equation}\label{eq1}
\dim(X)\geq \dim(F+R)=2m-1.
\end{equation}

\begin{prop}\label{2.3}
Let $X=\sum_{F\in\msf}F$.
\begin{itemize}
\item[\rm(i)]If $\dim(X)=2m-1$, then $|\msf|\leq \gc{2m-1}{m}$ and equality holds if and only if $\msf=\gc{X}{m}$.
\item[\rm(ii)]If $\dim(X)\geq 2m$, then
\begin{equation}\label{eq2.3.1}
|\msf|\leq \gc{m-1}{1}\gc{m}{1}^{m-1}+\gc{m-1}{1}^2\gc{2m-3}{m-2}.
\end{equation} 
\end{itemize}
\end{prop}
\proof(i) The fact that $\gc{X}{m}$ is an intersecting family with $\tau(\gc{X}{m})=m$ implies the desired result. 

(ii) We divide our discussion into two cases.
\medskip

\noindent {\bf Case 1}: $|\msf_T|\leq\gc{m-1}{1}$ for any $T\in\gc{V}{m-1}$.

Pick $F_0\in\msf$. Since $\dim(F\cap F_0)\geq 1$ for any $F\in\msf$, we have $\msf=\bigcup_{E\in\gc{F_0}{1}}\msf_E$. Corollary \ref{2.2} implies that
$$|\msf|\leq\sum_{E\in\gc{F_0}{1}}|\msf_E|\leq\gc{m-1}{1}\gc{m}{1}^{m-1}.$$
Hence, (\ref{eq2.3.1}) holds. 
\medskip

\noindent {\bf Case 2}: There exists a $T\in\gc{V}{m-1}$ with $|\msf_T|>\gc{m-1}{1}$.

Write $Y=\sum_{F\in\msf_T}F$, $\msa=\{A\in\msf:A\cap T=0\}$, $\msb=\{B\in\msf: B\not\subseteq Y\}$ and $\msc=\{C\in\msf:C\subseteq Y\}$.

We claim that $\dim(Y)=2m-1$ and $\msa\subseteq\gc{Y}{m}$. By Corollary 1.9 in \cite{WZX}, one gets
$$\gc{m-1}{1}<|\msf_T|\leq\left|\left\{F\in\gc{Y}{m}: T\subseteq F\right\}\right|=\gc{\dim(Y)-(m-1)}{1},$$
which implies that $\dim(Y)\geq 2m-1$. For any $F\in\msf_T$ and $A\in\msa$, it's routine to check that $F=T+(F\cap A)$ holds. Therefore, 
$$Y=\sum_{F\in\msf_T}F=\sum_{F\in\msf_T}\left(T+(F\cap A)\right)=T+\sum_{F\in\msf_T}(F\cap A)\subseteq T+A.$$
It follows that $\dim(Y)\leq \dim(T+A)=2m-1$, and so $Y=T+A$ holds for any $A\in\msa$. Hence, our claim is valid.

Pick $A_0\in\msa$. For any $B\in\msb$, the fact that $B\not\in\msa$ implies that $\dim(B\cap T)\geq 1$ and $\dim(B\cap A_0)\geq 1$. By Corollary \ref{2.2}, one gets
\begin{equation*}
|\msb|\leq\sum_{\scriptscriptstyle E_1\in\gc{A_0}{1},E_2\in\gc{T}{1}}|\msb_{E_1+E_2}|\leq\sum_{\scriptscriptstyle E_1\in\gc{A_0}{1},E_2\in\gc{T}{1}}|\msf_{E_1+E_2}|\leq \gc{m-1}{1}\gc{m}{1}^{m-1}.
\end{equation*}

Pick $B_1\in\msb$. Note that $1\leq \dim(B_1\cap Y)\leq m-1$. Since $\tau(\msf)=m$, there exists $B_2\in\msf$ such that $(B_1\cap Y)\cap B_2=0$, which implies that $B_2\in\msb$ and $1\leq\dim(B_2\cap Y)\leq m-1$. Observe that for any $C\in\msc$, $\dim(C\cap (B_i\cap Y))\geq 1$ for $i=1,2$. Since $(B_1\cap Y)\cap(B_2\cap Y)=0$ and $|\msc_D|\leq\gc{2m-3}{m-2}$ for any $D\in\gc{V}{2}$, we have
$$|\msc|\leq\sum_{\scriptscriptstyle{E_1\in\gc{B_1\cap Y}{1},E_2\in\gc{B_2\cap Y}{1}}}|\msc_{E_1+E_2}|\leq\gc{m-1}{1}^2\gc{2m-3}{m-2}.$$

Since $|\msf|=|\msb|+|\msc|$, (\ref{eq2.3.1}) holds.
$\qed$

{\bf{Proof of Theorem \ref{1.2}}}: By Proposition \ref{2.3}, it suffices to prove that 
\begin{equation}\label{eq2.3}
\gc{m-1}{1}\gc{m}{1}^{m-1}+\gc{m-1}{1}^2\gc{2m-3}{m-2}<\gc{2m-1}{m}.
\end{equation}

It is routine to check that (\ref{eq2.3}) holds for $m=3$. Suppose $m\geq 4$. 
Observe that
\begin{align}
&\gc{2m-1}{m}=\prod_{i=0}^{m-1}\frac{q^{2m-1-i}-1}{q^{m-i}-1}>q^{m(m-1)},\label{eq13}\\
&\gc{2m-1}{m}=\frac{q^{2m-1}-1}{q^m-1}\cdot\frac{q^{2m-2}-1}{q^{m-1}-1}\gc{2m-3}{m-2}>q^{2(m-1)}\gc{2m-3}{m-2}.\label{eq14}
\end{align}

By Bernoulli's inequality, one gets
$$\left(1-\frac{1}{q}\right)^{m-3}\geq 1-\frac{m-3}{q}>1-\frac{m-3}{q-1},$$
which implies that
\begin{equation}\label{eq12}
\frac{q^{m-1}}{(q-1)^{m-2}}<\frac{q^2}{q-m+2}\leq \frac{q^2}{q-2}\leq q(q-2).
\end{equation} 

Note that the left hand side of (\ref{eq2.3}) is less than
\begin{equation}\label{eq7}
\frac{q^{m-1}}{(q-1)^{m-2}}\cdot \frac{1}{(q-1)^2}\cdot q^{m(m-1)}+\frac{q^{2(m-1)}}{(q-1)^2}\gc{2m-3}{m-2}.
\end{equation}
Combining (\ref{eq13}), (\ref{eq14}) and (\ref{eq12}), we obtain that (\ref{eq7}) is less than $\displaystyle{{2m-1}\brack{m}}$, as desired.
$\qed$
\section{Singular linear spaces}
For fixed integers $n,l$ with $n>0$ and $l\geq 0$, let $V$ be an $(n+l)$-dimensional vector space over $\mathbb{F}_q$. A flag in $V$ is a sequence $(W_1,\ldots,W_r)$ of subspaces of $V$ such that
$$W_1\subseteq\cdots\subseteq W_r=V.$$
A \emph{parabolic subgroup} of general linear group $GL(V)$ of $V$ is the stabilizer of some flag in $V$. In particular, the parabolic subgroup $G$ determined by a flag $(W_1,W_2)$ with $\dim(W_1)=l$ is called a \emph{singular linear group}. The space $V$ together with the action of $G$ is called a \emph{singular linear space}. We say that an $m$-dimensional subspace $P$ is of type $(m,k)$ if $\dim(P\cap W_1)=k$. Denote the set of all subspaces of type $(m,k)$ in $V$ by $\msm(m,k;n+l,n)$. From Lemma 2.1 in \cite{KWW}, note that if $\msm(m,k;n+l,n)$ is non-empty, then it forms an orbit of subspaces under $G$.

The maximum size of an intersecting family $\msf\subseteq\msm(m,k;n+l,n)$ was determined in \cite{LO,Huang}. It was showed that an intersecting family with maximum size has covering number $1$. In \cite{GC}, we determined the maximum size of an intersecting family $\msf\subseteq\msm(n,0;n+l,n)$ whose covering number is more than $1$. In this section, we shall characterize the largest intersecting family $\msf\subseteq\msm(m,k;n+l,n)$ with $\tau(\msf)=m$. 

\begin{lem}\label{3.7}
Suppose $0\leq d\leq a-b$. Let $A\in\gc{V}{a}$ and $B\in\gc{V}{b}$ with $\dim(A\cap B)=c$. Then there exists $S\in\gc{A+B}{b-c+d}$ such that $\dim(S\cap A)=d$ and $S\cap B=0$.
\end{lem}
\proof
Denote $C=A\cap B$. If $B\subseteq A$, then the result is directed. Suppose $B\not\subseteq A$. Pick $A_1\in\gc{A}{b-c+d}$ and $B_1\in\gc{B}{b-c}$ such that $A_1\cap C=B_1\cap C=0$. Let $\alpha_1,\ldots,\alpha_{b-c+d}$ be a base of $A_1$ and $\beta_1,\ldots, \beta_{b-c}$ be a base of $B_1$. Write
\begin{equation*}
S=
\begin{cases}
\langle \alpha_1+\beta_1,\ldots, \alpha_{b-c}+\beta_{b-c},\alpha_{b-c+1},\ldots,\alpha_{b-c+d}\rangle,&~\text{if}~d\geq 1,\\
\langle \alpha_1+\beta_1,\ldots, \alpha_{b-c}+\beta_{b-c}\rangle,&~\text{if}~d=0.
\end{cases}
\end{equation*}
It's routine to check that $S$ is the desired subspace.
$\qed$

Let $N(m_1,k_1;m,k;n+l,n)$ be the number of subspaces of type $(m_1,k_1)$ contained in a given subspace of type $(m,k)$ in $V$.
 
\begin{lem}{\rm{(\cite[Lemma 2.2]{KW})}}\label{3.2}
$N(m_1,k_1;m,k;n+l,n)>0$ if and only if 
\begin{equation}\label{eq3.2}
0\leq k_1\leq k\leq l,~0\leq m_1-k_1\leq m-k\leq n.
\end{equation}
Moreover, if $N(m_1,k_1;m,k;n+l,n)>0$, then
\begin{equation*}
N(m_1,k_1;m,k;n+l,n)=q^{(m_1-k_1)(k-k_1)}\gc{m-k}{m_1-k_1}\gc{k}{k_1}.
\end{equation*}
\end{lem}

In the remaining of this paper, suppose that $\msf\subseteq\msm(m,k;n+l,n)$ is an intersecting family with $\tau(\msf)=m$. Write $X=\sum_{F\in\msf}F$. By (\ref{eq1}), one gets $\dim(X)\geq 2m-1$. Let $\gc{X}{m,k}$ be the family of all subspaces of type $(m,k)$ contained in $X$. 

\begin{lem}\label{3.3}
Suppose $\dim(X)=2m-1$. Then $|\msf|\leq N(m,k;2m-1,t;n+l,n)$ and equality holds if and only if $\msf=\gc{X}{m,k}$, where $X$ is of type $(2m-1,t)$ with 
\begin{equation}\label{eq99}
t=
\begin{cases}
\max\{0,2m-1-n\},~&\text{if}~k=0,\\ 
m+k-1,~&\text{if}~k\geq 1.
\end{cases}
\end{equation}
\end{lem}
\proof
Write $X'=X\cap W_1$, $\dim(X')=t$. Then
\begin{equation}\label{eq4}
|\msf|\leq\left|\gc{X}{m,k}\right|=N(m,k;2m-1,t;n+l,n).
\end{equation}

If $k=0$, then by Lemma \ref{3.2}, one gets that $N(m,0;2m-1,t;n+l,n)$ reach the maximum value exactly at $t=\max\{0,2m-1-n\}$. Suppose $k\geq 1$. Observe that $k\leq t\leq m+k-1$ by Lemma \ref{3.2}. Pick $Z\in\gc{X'}{t-k+1}$. Then $Z$ intersects every $k$-subspace of $X'$. Note that $\dim(F\cap X')=k$ holds for any $F\in\msf$. Thus, $Z$ intersects every element of $\msf$, which implies that $m=\tau(\msf)\leq\dim(Z)=t-k+1$,
 and so $t=m+k-1$. Hence, we characterize the exactly case when $|\gc{X}{m,k}|$ reach the maximum value.

In order to finish our proof, it suffices to show that $\tau(\gc{X}{m,k})=m$. For any $I\in\gc{X}{i}$ with $1\leq i\leq m-1$, write $H=I\cap X$, $\dim(H)=h$ and $\dim(X'\cap H)=j$. Then $\dim(X'+H)=t+h-j$. We only need to prove that there exists an $F\in\gc{X}{m,k}$ such that $F\cap I=0$. 

\noindent{\bf Case 1}: $k=0$.

Since $m\leq n$, we have $t\leq m-1$. Note that

\begin{equation}\label{eq100}
\left.
\begin{aligned}
&\max\{t-j,h-j\}<m,\\
&\max\{m-(t-j),m-(h-j)\}\leq 2m-1-(t+h-j).
\end{aligned}
\right\}
\end{equation}

\noindent{\bf Case 1.1}: $t\leq h$. 

By Lemma \ref{3.7}, there exists a $D_1\in\gc{X'+H}{t-j}$ such that $D_1\cap X'=D_1\cap H=0$. By (\ref{eq100}), there exists an $F\in\gc{X}{m}$ such that $F\cap (X'+H)=D_1$, which implies that
\begin{align*}
&F\cap W_1=F\cap X'=D_1\cap X'=0,\\
&F\cap I=F\cap H=D_1\cap H=0.
\end{align*}
Hence, $F$ is the desired subspace. 

\noindent{\bf Case 1.2}: $t>h$. 

Similar to Case 1.1, we may find a desired subspace $F$.

\noindent{\bf Case 2}: $k\geq 1$.

Observe that $t=m+k-1\geq h+k$. By Lemma \ref{3.7}, there exists an $S\in\gc{X'+H}{h+k-j}$ such that $S\cap H=0$ and $\dim(S\cap X')=k$. Since $(X'+H)\subseteq X$, one gets that $m\geq h+k-j$, which implies that there exists an $F\in\gc{V}{m}$ such that $F\cap(X'+H)=S$. Note that
\begin{align*}
&F\cap W_1=F\cap X'=S\cap X',\\
&F\cap I=F\cap H=S\cap H=0.
\end{align*}
Hence, $F$ is the desired subspace. 
$\qed$

\begin{thm}\label{3.1}
Suppose $q\geq m+2\geq 5$. Let $\msf\subseteq\msm(m,k;n+l,n)$ be an intersecting family with $\tau(\msf)=m$. Then $|\msf|\leq N(m,k;2m-1,t;n+l,n)$ and equality holds if and only if $\msf=\gc{X}{m,k}$, where $X$ is of type $(2m-1,t)$, and $t$ is as in (\ref{eq99}).
\end{thm} 
\proof
By (\ref{eq2.3.1}) and Lemma \ref{3.3}, it suffices to prove that
\begin{equation}\label{eq233}
\gc{m-1}{1}\gc{m}{1}^{m-1}+\gc{m-1}{1}^2\gc{2m-3}{m-2}<N(m,k;2m-1,t;n+l,n).
\end{equation}

Suppose $a>b\geq 1$. Observe that the function $\dfrac{q^{a-x}-1}{q^{b-x}-1}$ monotonically increases in the interval $[0,b-1]$, which implies that
$$q^{a-b}<\frac{q^{a-i}-1}{q^{b-i}-1}\leq\gc{a-b+1}{1},~i\in\{0,\ldots,b-1\}.$$
Therefore, we obtain that
\begin{equation}\label{eq10}
q^{b(a-b)}<\gc{a}{b}\leq \gc{a-b+1}{1}^{b}.
\end{equation}

By the second inequality in (\ref{eq10}), we have $\gc{2m-3}{m-2}\leq\gc{m}{1}^{m-2}$, which implies that the left hand side of (\ref{eq233}) is less than
\begin{equation}\label{eq11}
\frac{q^{m-1}}{(q-1)^m}\cdot q^{m(m-1)}+\frac{q^{m-2}}{(q-1)^m}\cdot q^{m(m-1)}.
\end{equation}
Since 
$$(q-1)^m>q^m-m q^{m-1}\geq q^m-(q-2)q^{m-1}>q^{m-1}+q^{m-2}$$
by Bernoulli's inequality, (\ref{eq11}) is less than $q^{m(m-1)}$. On the other hand, by Lemma \ref{3.2} and the first inequality in (\ref{eq10}), 
$$q^{m(m-1)}\leq N(m,k;2m-1,t;n+l,n).$$
Hence, (\ref{eq233}) holds.
$\qed$
\section*{Acknowledgement}

This research is supported by NSFC (11671043).

\end{document}